\documentclass{elsart}

\usepackage{amsmath,amssymb}
\usepackage{natbib}
\usepackage{tikz}

\newcommand{\F}{\mathbb{F}}

\newcommand{\ga}{\alpha}

\renewcommand{\phi}{\varphi}

\newcommand{\set}[1]{\{#1\}}

\newcommand{\ev}{\mathrm{ev}}
\newcommand{\wt}{\mathrm{wt}}
\newcommand{\RS}{\mathrm{RS}}
\newcommand{\LT}{\mathrm{lt}}
\newcommand{\LC}{\mathrm{lc}}
\newcommand{\ydeg}{\text{$y$-$\deg$}}
\newcommand{\wdeg}{\deg_u}

\DeclareMathOperator{\mult}{mult}

\begin{document}

\begin{frontmatter}

\title{List Decoding of Reed-Solomon Codes \\ from a Gr\"obner Basis Perspective}

\author{Kwankyu~Lee\corauthref{cor}}
\ead{kwankyu@kias.re.kr}
\address{School of Computational Sciences \\ Korea Institute for Advanced Study, Seoul, Korea}
\author{Michael E.~O'Sullivan}
\ead{mosulliv@math.sdsu.edu}
\corauth[cor]{Corresponding author.}
\address{Department of Mathematics and Statistics \\ San Diego State University, San Diego, USA}

\begin{abstract}
The interpolation step of Guruswami and Sudan's list decoding of Reed-Solomon codes poses the problem of finding the minimal polynomial of an ideal with respect to a certain monomial order. An efficient algorithm that solves the problem is presented based on the theory of Gr\"obner bases of modules. In a special case, this algorithm reduces to a simple Berlekamp-Massey-like decoding algorithm. 
\end{abstract}


\end{frontmatter}

\section{Introduction}

Interpreting the key equation of \citet{welch1986} as a problem of finding a plane curve interpolating points with a certain weighted degree constraint, \citet{sudan1997} developed list decoding of Reed-Solomon codes. Soon afterward, using the concept of multiplicity at a point on an algebraic curve, \citet{guruswami1999} improved Sudan's list decoding so that it is capable of correcting more errors than conventional decoding algorithms for all rates of Reed-Solomon codes. The list decoding consists of two steps: the interpolation step and the root-finding step, each of which poses a problem that can be solved in various ways. Since the interpolation problem can be solved by finding a solution of a system of linear equations over a field, they simply asserted the existence of a polynomial time algorithm solving the interpolation problem, thus leaving it as an open problem to search for an efficient interpolation algorithm. 
  
Several authors, including \citet{nielsen2000}, \citet{fitzpatrick2002}, \citet{mceliece2003} in his presentation of K\"otter's algorithm, \citet{aleknovich2005}, and \citet{farr2005}, formulated the interpolation problem as a problem of finding the minimal polynomial, with respect to a weighted monomial order, of the ideal of polynomials interpolating certain points. Their interpolation algorithms, except Alekhnovich's, take basically a ``\emph{point by point}'' approach in the sense that they build a Gr\"obner basis of the ideal for points $\set{P_1,P_2,\dots,P_n}$ by recursively computing a Gr\"obner basis of the ideal for $\set{P_1,\dots,P_i}$ while $i$ increases from $1$ to $n$. In this paper, we also take the Gr\"obner basis perspective, but employ a different strategy. We start with a set of generators of the module induced from the ideal for  $\set{P_1,P_2,\dots,P_n}$ and convert the generators to a Gr\"obner basis of the module, in which the minimal polynomial is found. This results in an efficient algorithm solving the interpolation problem.

In Section 2, we briefly review Guruswami and Sudan's list decoding of Reed-Solomon codes. A more detailed treatment can be found in \citet{mceliece2003} and in \citet{guruswami2005}. In Section 3, we formulate the interpolation problem in a Gr\"obner basis perspective. The basics of Gr\"obner bases that we assumed in this paper can be found in \citet{cox1997,cox2005}. In Sections 4--6, our interpolation algorithm is presented and analyzed. In Section 7, we treat the special case of multiplicity one and list size one. 

\section{List Decoding of Reed-Solomon Codes}

Let $\F$ be a finite field. Denote by $\F[x]_{s}$ the set of polynomials with degree $<s$, which is an $s$-dimensional subspace of $\F[x]$ as $\F$-vector spaces. Fix $n$ distinct points $\ga_1,\ga_2,\dots,\ga_n$ from $\F$. Note that the evaluation map $\ev:\F[x]_n\to\F^n$ defined by $f\mapsto(f(\ga_1),f(\ga_2),\dots,f(\ga_n))$ is an isomorphism of $\F$-vector spaces. The inverse map $\ev^{-1}$ is given by Lagrange interpolation as follows. Define
\begin{equation}\label{equcss}
	\tilde h_i=\prod_{j=1,j\neq i}^n(x-\ga_j),\text{ and } h_i=\tilde h_i(\ga_i)^{-1}\tilde h_i
\end{equation}
so that $h_i(\ga_j)=1$ if $j=i$, and $0$ otherwise. For $v=(v_1,v_2,\dots,v_n)\in\F^n$, we write
\[
	h_v=\sum_{i=1}^nv_ih_i\in\F[x]_n
\]
so that $h_v=\ev^{-1}(v)$. For $k<n$, the Reed-Solomon code $\RS(n,k)$, defined as the image of $\F[x]_k$ by $\ev$, is an $[n,k]$ linear code over $\F$. 

For $f\in\F[x,y]$ and $u\ge 1$, denote by $\deg_u(f)$ the $(1,u)$-weighted degree of $f$. That is, variables $x$ and $y$ are assigned weights $1$ and $u$ respectively, and for a monomial $x^iy^j$, we define $\deg_u(x^iy^j)=i+uj$. For a polynomial $f$, we define $\deg_u(f)$ as the maximal $\deg_u(x^iy^j)$ for monomials $x^iy^j$ occurring in $f$.

A nonzero polynomial in $\F[x,y]$ defines a curve on the plane $\F^2$. The multiplicity of a curve $f$ at the origin is defined to be the smallest $m$ such that a monomial of total degree $m$ occurs in the polynomial $f$. The multiplicity of a curve $f$ at an arbitrary point $P=(a,b)$ is defined as the multiplicity of the curve $f_P$ at the origin, where $f_P=f(x+a,y+b)$, and denoted by $\mult_P(f)$.

Suppose that some codeword of $\RS(n,k)$ was sent through a noisy channel, and the vector $v\in\F^n$ is received by hard-decision on the channel ouput. 
For each $1\le i\le n$, let $P_i$ denote the point $(\ga_i,v_i)$ on the plane $\F^2$. Now for $m\ge 1$, define 
\[
	I_{v,m}=\set{f\in\F[x,y]\mid\text{$\mult_{P_i}(f)\ge m$ for $1\le i\le n$}}\cup\set{0},
\]
which is an ideal of $\F[x,y]$. \cite{guruswami1999} proved
\begin{prop}\label{propdjqq}
Let $v\in\F^n$ be the received vector. Suppose that $f\in I_{v,m}$ is nonzero. Let $w=\deg_{k-1}(f)$. If $c$ is a codeword of $\RS(n,k)$ satisfying
\[
	\wt(v-c)<n-\frac{w}{m},
\]
then $h_c$ is a root of $f$ as a polynomial in $y$ over $\F[x]$.
\end{prop}
This proposition is the basis of their list decoding algorithm. We recall that the goal of the interpolation step of list decoding is to find a polynomial in $I_{v,m}$ having the smallest $(1,k-1)$-weighted degree. Having the same weighted degree, the one with smaller degree in $y$ is preferred because this reduces the work of the root-finding algorithm.

\section{Gr\" obner Basis Perspective}

We observe that if $I$ is an ideal of $\F[x,y]$, then the minimal polynomial of $I$ with respect to a monomial order $>$ is the minimal element of any Gr\"obner basis of the ideal $I$ with respect to $>$. This is a direct consequence of the definition of Gr\"obner bases. Let $Q$ be the minimal polynomial of $I_{v,m}$ with respect to the monomial order $>_{k-1}$ of $\F[x,y]$. As we observed, we can find $Q$ by computing a Gr\" obner basis of $I_{v,m}$ with respect to $>_{k-1}$. However, computing a Gr\"obner basis of an ideal is generally a task of high complexity. We overcome this difficulty by using the theory of Gr\"obner bases of modules.

Let $l$ be a positive integer. Let $\F[x,y]_l=\set{f\in\F[x,y]\mid\ydeg(f)\le l}$. We view $\F[x,y]_l$ as a free module over $\F[x]$ with a free basis $1,y,y^2,\dots,y^l$. Monomials of the module $\F[x,y]_l$ consist of $x^iy^j$ with $i\ge 0$ and $0\le j\le l$.  

Note that a monomial order $>$ on the ring $\F[x,y]$ naturally induces a monomial order on the module $\F[x,y]_l$, which we also denote by $>$. The notions of $(1,u)$-weighted degrees and $y$-degrees of monomials or polynomials in $\F[x,y]$ carry over to $\F[x,y]_l$. Thus $>_u$ is a monomial order on the module $\F[x,y]_l$. The notion of the minimal polynomial of a submodule of $\F[x,y]_l$ is defined in the same way as for an ideal of $\F[x,y]$. 

For $l\ge 1$, we define 
\[
	I_{v,m,l}=I_{v,m}\cap\F[x,y]_l.
\]
Then $I_{v,m,l}$ is a submodule of $\F[x,y]_l$. The minimal polynomial $Q$ of $I_{v,m}$ with respect to $>_{k-1}$ is also the minimal polynomial of $I_{v,m,l}$ with respect to $>_{k-1}$ if $l$ is as large as the $y$-degree of $Q$. This enables us to find $Q$ by computing a Gr\" obner basis of the submodule $I_{v,m,l}$ of the free module $\F[x,y]_l$ with respect to $>_{k-1}$. This task turns out to be much easier than that of computing a Gr\"obner basis of the ideal $I_{v,m}$ because there is a simple criterion of Gr\"obner bases for a submodule of $\F[x,y]_l$. The following is a trivial application of Buchberger's $S$-pair criterion.

\begin{prop}\label{propawjw}
Let $S$ be a submodule of $\F[x,y]_l$ with a monomial order $>$. Suppose that $\set{g_0,g_1,\dots,g_s}$ generates $S$. If $y$-degrees of leading terms of $g_i$ for $0\le i\le s$ are all distinct, then $\set{g_0,g_1,\dots,g_s}$ is a Gr\"obner basis of $S$ with respect to $>$.
\end{prop}

It is easy to identify a set of generators of $I_{v,m,l}$, from which we compute a Gr\"obner basis. First we present a natural set of generators for the ideal $I_{v,m}$.

\begin{prop}
As an ideal of $\F[x,y]$,
\[
	I_{v,m}=\langle y-h_v,\eta\rangle^m
	 =\langle(y-h_v)^i\eta^{m-i}\mid 0\le i\le m\rangle,
\]
where $\eta=\prod_{j=1}^n(x-\ga_j)$.
\end{prop}

\begin{pf}
Let 
\[
	J=\langle(y-h_v)^i\eta^{m-i}\mid 0\le i\le m\rangle.
\]
Then $J\subset I_{v,m}$ since each generator of $J$ is clearly in $I_{v,m}$. To prove the reverse relation, let $f\in I_{v,m}$. By division with respect to $y$, we write
\[
	f=g_0(y-h_v)^m+f_0
\]
with $\ydeg(f_0)<m$. Note that $f_0\in I_{v,m}$ since $J\subset I_{v,m}$. Let $d=\ydeg(f_0)$, and write
\[
	f_0=g(y-h_v)^d+f_1
\]
with $g\in\F[x]$ and $\ydeg(f_1)<d$. Observe that for $1\le s\le n$, 
\[
\begin{split}
	f_0(x+\ga_s,y+v_s)&=g(x+\ga_s)(y+v_s-h_v(x+\ga_s))^d+f_1(x+\ga_s,y+v_s) \\
		&=g(x+\ga_s)(y^d+\cdots)+f_1(x+\ga_s,y+v_s).
\end{split}
\]
As $f_0$ has multiplicity at least $m$ at $P_s$, $f_0(x+\ga_s,y+v_s)$ has no monomial of total degree less than $m$. Because $\ydeg(f_1)<d$, we see that $g(x+\ga_s)$ must be divisible by $x^{m-d}$, which implies $(x-\ga_s)^{m-d}$ divides $g(x)$. Therefore we can write with some $g_1\in\F[x]$,
\[
	f_0=g_1(y-h_v)^d\prod_{j=1}^n(x-\ga_j)^{m-d}+f_1.
\]
We continue this until we eventually have $f_i=0$ as $y$-degrees are decreasing. Then $f\in J$. Hence $I_{v,m}=J$.
\end{pf}
 
\begin{cor}\label{propskqq}
Let $l\ge m$. As a submodule of $\F[x,y]_l$ over $\F[x]$,
\[
	I_{v,m,l}=\langle (y-h_v)^i\eta^{m-i},y^{i'-m}(y-h_v)^m
	\mid0\le i\le m,m< i'\le l\rangle,
\]
where $\eta=\prod_{j=1}^n(x-\ga_j)$.
\end{cor}

We need an upper bound on $\ydeg(Q)$ to set $l$. As in \cite{guruswami1999}, using the following 
\begin{prop}\label{propsjqx}
Let $S$ be a subset of exponents of monomials of $\F[x,y]$. If $|S|$ is at least
\begin{equation}\label{equxqwn}
	N=n\binom{m+1}{2}+1,
\end{equation}
then there is a set of coefficients $f_{ij}$ such that $f=\sum_{(i,j)\in S}f_{ij}x^iy^j$ is a nonzero polynomial in $I_{v,m}$.
\end{prop}
and counting the monomials ordered in $>_{k-1}$, we can get
\begin{equation}\label{dmxcxn}
	\ydeg(Q)<\sqrt{\frac{2N}{k-1}+\frac{1}{4}}-\frac{1}{2}.
\end{equation}
We refer to the original source for a proof of the proposition and a detailed derivation of the upper bound. The derivation also implies that the upper bound is larger than $m$ as $N$ is not less than the number of monomials from $1$ to $y^m$ inclusive ordered with respect to $>_{k-1}$. Henceforth we let $l$ be the largest integer less than the upper bound.

In the next section, we present an algorithm converting the set of generators of $I_{v,m,l}$ given in Corollary \ref{propskqq} to a Gr\"obner basis with respect to $>_{k-1}$. 

\section{A Gr\"obner Basis Algorithm}

Let $S$ be a submodule of $\F[x,y]_l$ over $\F[x]$. Suppose that $\set{g_0,g_1,\dots,g_l}$
is a set of generators of $S$ and satisfy $\ydeg(g_i)=i$ for $0\le i\le l$ (equivalently the generators form a Gr\"obner basis with respect to the lexicographical order with $y>x$). Fix a monomial order $>_u$ on $\F[x,y]_l$. The following algorithm computes a Gr\" obner basis of $S$ with respect to $>_u$ from $g_0,g_1,\dots,g_l$.

\textbf{Algorithm G.} Let $g_i=\sum_{j=0}^la_{ij}y^j$ for $0\le i\le l$ during the execution of the algorithm.

\begin{enumerate}
\item[G1.] Set $r\leftarrow 0$.
\item[G2.] Increase $r$ by $1$. If $r\le l$, then proceed; otherwise go to step G6.
\item[G3.] Find $s=\ydeg(\LT(g_r))$. If $s=r$, then go to step G2.
\item[G4.] Set $d\leftarrow\deg(a_{rs})-\deg(a_{ss})$ and  $c\leftarrow\LC(a_{rs})\LC(a_{ss})^{-1}$.
\item[G5.] (a) If $d\ge 0$, then set 
\[
	g_r\leftarrow g_r-cx^dg_s.
\]
(b) If $d<0$, then set, storing $g_s$ in a temporary variable,
\[
	g_s\leftarrow g_r,\quad g_r\leftarrow x^{-d}g_r-cg_s.
\]
Go back to step G3.
\item[G6.] Output $\set{g_0,g_1,\dots,g_l}$ and the algorithm terminates.  
\end{enumerate}

The goal of the algorithm is to inductively process $g_0,g_2,\dots,g_l$ such that they still generate $S$ and $\ydeg(\LT(g_i))=i$ for $0\le i\le l$, so that $\set{g_0,g_1,\dots,g_l}$ is a Gr\"obner basis of $S$ by Proposition \ref{propawjw}. Note that initially we have
\[
	\left\{\begin{array}{cr}
	g_0=&a_{00}\\
	g_1=&a_{11}y+a_{10}\\
	g_2=&a_{22}y^2+a_{21}y+a_{20}\\
	\vdots\\
	g_l=&a_{ll}y^l+\cdots+a_{l2}y^2+a_{l1}y+a_{l0}
	\end{array}\right.
\]
After increasing $r$ by one in step G2, $\ydeg(\LT(g_i))=i$ for $0\le i\le r-1$. Then the algorithm processes $g_0,g_1,\dots,g_r$ by iterating steps G3--G5, until $\ydeg(\LT(g_i))=i$ for $0\le i\le r$. Observe that $g_s$ and $g_r$ are updated in a way that the new $g_0,g_1,\dots,g_l$ still generate the module $S$. When the algorithm terminates, we have $\ydeg(\LT(g_i))=i$ for $0\le i\le l$ as desired.

We may view Algorithm G as an optimized version of Buchberger's algorithm. However, to prove directly that the algorithm terminates and hence output a Gr\"obner basis is even easier.

\begin{prop}\label{propcjdc}
Fix $r>0$ and suppose $\ydeg(\LT(g_i))=i$ for $0\le i<r$. After a finite number of iterations through steps \textup{G3--G5}, it eventually happens that $\ydeg(\LT(g_i))=i$ for $0\le i\le r$. 
\end{prop}

\begin{pf}
Observe that the update (a) does not change the weighted degrees of $\LT(g_i)$ for $0\le i<r$ while the update (b) strictly decreases the weighted degree of $\LT(g_s)$ but keeps the weighted degrees of $\LT(g_i)$ for $0\le i<r$ with $i\neq s$. Therefore the update (b) could not occur infinitely many times. So from a certain point on, only the update (a) occurs. Now observe that the update (a) either strictly decreases the weighted degree of $g_r$ or otherwise the $y$-degree of $\LT(g_r)$ strictly decreases. Therefore the update (a) could not happen infinitely. Hence iterations must stop either by $g_r$ vanishing to zero or by the $y$-degree of $\LT(g_r)$ being $r$. However, the first case is not possible because $g_0,g_1,\dots,g_r$ form a rank $r+1$ free module over $\F[x]$.
\end{pf}

Unfortunately, the above proof does not allow us to estimate the complexity of the algorithm because we cannot know how many iterations of steps G3--G5 occur before the algorithm terminates. Hence we need to understand the behavior of the algorithm more carefully. 

\begin{prop}
Let $g_i = \sum_{j=0}^r a_{ij}y^j$ and $g_i'=\sum_{j=0}^r a_{ij}' y^j$, $0\le i\le r$ be the states of the algorithm before and after step \textup{G5}, respectively. Then for any non-identity permutation $\pi=(\pi_0,\pi_1,\dots,\pi_r)$,
\begin{align}\label{cjwqwdd}
	\sum_{i=0}^r \deg (a_{ii}') &> \sum_{i=0}^r \deg (a_{i\pi_i}'). 
\end{align}
Moreover if $d\ge 0$, then 
\begin{equation}\label{cjswdf}
 \deg(a_{rr}')=\deg(a_{rr})\text{ and }\deg_u(a_{rj}'y^j)\le\deg_u(a_{rs}y^s)
\end{equation}
for $j\le r$ with strict inequality for $j\ge s$. Similarly if $d<0$, then
\begin{equation}\label{cjswdg}
 \deg(a_{rr}')=\deg(a_{rr})-d\text{ and }\deg_u(a_{rj}'y^j)\le\deg_u(a_{rs}y^s)-d
\end{equation}
for $j\le r$ with strict inequality for $j\ge s$.
\end{prop}

\begin{pf}
For induction, let us assume 
\begin{align}\label{cjwqwde}
	\sum_{i=0}^r \deg (a_{ii}) &> \sum_{i=0}^r \deg (a_{i\pi_i}) 
\end{align}
for any non-identity permutation $\pi$. First consider the case $d=\deg(a_{rs})-\deg(a_{ss})\ge 0$, where $a_{rj}' = a_{rj} - cx^da_{sj}$. Applying \eqref{cjwqwde} to the transposition of $s$ and $r$, we get $\deg(a_{rr})+\deg(a_{ss})>\deg(a_{rs})+\deg(a_{sr})$, which implies the equality part of \eqref{cjswdf}. The inequality part of \eqref{cjswdf} follows by noting $\deg_u(a_{rj}y^j)\le\deg_u(a_{rs}y^s)$ and
\[
	\wdeg(x^da_{sj}y^j)=\wdeg(a_{rs}y^s)+(\wdeg(a_{sj}y^j)-\wdeg(a_{ss}y^s))\le\wdeg(a_{rs}y^s)
\]
for $j\le r$ with strict inequality for $j>s$ and by noting for $j=s$ the way that $c$ and $d$ is chosen. 

We turn to \eqref{cjwqwdd}. By what we proved, the left hand side of \eqref{cjwqwdd} equals $\sum_{i=0}^r \deg(a_{ii})$. For the right hand side, note that $\deg(a_{r\pi_r}')\le\deg(a_{r\pi_r})$ or $\deg(a_{r\pi_r}')=\deg(a_{rs})-\deg(a_{ss})+\deg(a_{s\pi_r})$. If the first case holds, then $\sum_{i=0}^r \deg (a_{i\pi_i}')\le\sum_{i=0}^r \deg (a_{i\pi_i})$, and \eqref{cjwqwdd} follows from \eqref{cjwqwde}. Supposing the second case, let $D_{ij}$ denote $\deg_u(a_{ij}y^j)$. Then \eqref{cjwqwdd} is equivalent to
\begin{equation}\label{eqppwk}
	\sum_{i}D_{ii}>\sum_{i\neq s,r}D_{i\pi_i}+D_{s\pi_s}+D_{s\pi_r}+D_{rs}-D_{ss}.
\end{equation}
To show \eqref{eqppwk}, we need to treat two cases depending on whether $s$ and $\pi_r$ are in the same orbit or not, with respect to the permutation $\pi$.	First suppose $s$ and $\pi_r$ are in the same orbit so that
\[
s\longrightarrow\pi_s\longrightarrow\cdots\longrightarrow\pi_r\longrightarrow\pi(\pi_r)\longrightarrow\cdots\longrightarrow\pi^{-1}(s)\longrightarrow s.
\]
Let $S=\set{\pi_r,\pi(\pi_r),\dots,\pi^{-1}(s)}$. Note that $S$ is empty if $\pi_r=s$. Now the right hand side of \eqref{eqppwk} equals
\begin{equation}\label{equcnxx}
\begin{split}
	&\sum_{i\in S}D_{i\pi_i}+\sum_{i\not\in S,\, i\neq s,r}\!\!\!D_{i\pi_i}+D_{s\pi_s}+D_{rs}-D_{ss}+D_{s\pi_r}\\
	&\quad\le\sum_{i\in S}D_{ii}+\sum_{i\not\in S,\, i\neq s,r}\!\!\!D_{i\pi_i}+D_{s\pi_s}+D_{rs}<\sum_{i}D_{ii}.
\end{split}
\end{equation}
Here the first inequality holds since $D_{ii}\ge D_{ij}$ for $0\le i\le r-1$ and $0\le j\le r$ by the algorithm. The second strict inequality follows from \eqref{cjwqwde} as we can check that the second indices of the terms in the middle expression of \eqref{equcnxx} are all distinct by the definition of $S$. 

If $s$ and $\pi_r$ are not in the same orbit, then we have
\[
s\longrightarrow\pi_s\longrightarrow\pi(\pi_s)\longrightarrow\cdots\longrightarrow\pi^{-1}(s)\longrightarrow s,
\]
and let $S=\set{\pi_s,\pi(\pi_s),\dots,\pi^{-1}(s)}$. Note that $S$ is empty if $\pi_s=s$. Now the right hand side of \eqref{eqppwk} equals
\[
\begin{split}
	&\sum_{i\in S}D_{i\pi_i}+\sum_{i\not\in S,\, i\neq s,r}\!\!\!D_{i\pi_i}+D_{s\pi_r}+D_{rs}-D_{ss}+D_{s\pi_s}\\
	&\quad\le\sum_{i\in S}D_{ii}+\sum_{i\not\in S,\, i\neq s,r}\!\!\!D_{i\pi_i}+D_{s\pi_r}+D_{rs}<\sum_{i}D_{ii},
\end{split}
\]
where the inequalities are justified by similar arguments as above.

Let us now consider the case $d<0$, where $a_{sj}'=a_{rj}$ and $a_{rj}'=x^{-d}a_{rj}-ca_{sj}$. We can verify \eqref{cjswdg} in a similar way to the case $d\ge 0$, so we turn to \eqref{cjwqwdd}. Since $\deg(a_{rr}')=\deg(a_{ss})-\deg(a_{rs})+\deg(a_{rr})$, the left hand side of \eqref{cjwqwdd} equals $\sum_{i=0}^r\deg(a_{ii})$. As $a_{r\pi_r}'=x^{-d}a_{r\pi_r}-ca_{s\pi_r}$, we have $\deg(a_{r\pi_r}')\le\deg(a_{s\pi_r})$ or $\deg(a_{r\pi_r}')=\deg(a_{ss})-\deg(a_{rs})+\deg(a_{r\pi_r})$. If the first case holds, then $\sum_{i=0}^r\deg(a_{i\pi_i}')\le\sum_{i\neq s,r}\deg(a_{i\pi_i})+\deg(a_{r\pi_s})+\deg(a_{s\pi_r})$, and \eqref{cjwqwdd} follows from \eqref{cjwqwde}. Suppose the second case, and let $D_{ij}$ denote $\deg_u(a_{ij}y^j)$. Note that \eqref{cjwqwdd} is equivalent to 
\begin{equation}\label{equchqw}
	\sum_{i}D_{ii}>\sum_{i\neq s,r}D_{i\pi_i}+D_{r\pi_s}+D_{r\pi_r}+D_{ss}-D_{rs}.
\end{equation}
To show this, we treat two cases depending on whether $s$ and $\pi_r$ are in the same orbit or not, with respect to the permutation $\pi$.	First suppose $s$ and $\pi_r$ are in the same orbit so that
\[
s\longrightarrow\pi_s\longrightarrow\cdots\longrightarrow\pi_r\longrightarrow\pi(\pi_r)\longrightarrow\cdots\longrightarrow\pi^{-1}(s)\longrightarrow s.
\]
Let $S=\set{\pi_r,\pi(\pi_r),\dots,\pi^{-1}(s)}$. Note that $S$ is empty if $\pi_r=s$. Now the right hand side of \eqref{equchqw} equals
\begin{equation}\label{equcnvf}
\begin{split}
	&\sum_{i\in S}D_{i\pi_i}+\sum_{i\not\in S,\, i\neq s,r}\!\!\!D_{i\pi_i}+D_{r\pi_s}+D_{ss}-D_{rs}+D_{r\pi_r}\\
	&\quad\le\sum_{i\in S}D_{ii}+\sum_{i\not\in S,\, i\neq s,r}\!\!\!D_{i\pi_i}+D_{r\pi_s}+D_{ss}<\sum_{i}D_{ii}.
\end{split}
\end{equation}
The first inequality holds since $D_{ii}\ge D_{ij}$ for $0\le i\le r-1$ and $0\le j\le r$ and $D_{rs}\ge D_{rj}$ for $0\le j\le r$ by the way in which $s$ is chosen. The second strict inequality follows from \eqref{cjwqwde} since we can check that the right indices of terms in the middle expression of \eqref{equcnvf} are all distinct. 

If $s$ and $\pi_r$ are not in the same orbit, then we have
\[
\pi_s\longrightarrow\pi(\pi_s)\longrightarrow\cdots\longrightarrow\pi^{-1}(s)\longrightarrow s,
\]
and let $S=\set{\pi_s,\pi(\pi_s),\dots,\pi^{-1}(s)}$. Note that $S$ is empty if $\pi_s=s$. Now the right hand side of \eqref{equchqw} equals
\[
\begin{split}
	&\sum_{i\in S}D_{i\pi_i}+\sum_{i\not\in S,\, i\neq s,r}\!\!\!D_{i\pi_i}+D_{r\pi_r}+D_{ss}-D_{rs}+D_{r\pi_s}\\
	&\quad\le\sum_{i\in S}D_{ii}+\sum_{i\not\in S,\, i\neq s,r}\!\!\!D_{i\pi_i}+D_{r\pi_r}+D_{ss}<\sum_{i}D_{ii},
\end{split}
\]
where the inequalities hold by the same reasons as above.
\end{pf}

\begin{cor}\label{fdhwkd}
With the notation of the proposition, we have
\[ 
	\deg_u(\LT(g_r'))-\deg_u(a_{rr}'y^r)\le\deg_u(\LT(g_r))-\deg_u(a_{rr}y^r)
\]
If equality holds, then $\ydeg(\LT(g_r'))<\ydeg(\LT(g_r))$.
\end{cor}

\begin{pf}
Note that $\deg_u(\LT(g_r))=\deg_u(a_{rs}y^s)$. Then the assertions are immediate from \eqref{cjswdf} and \eqref{cjswdg}.
\end{pf}

\section{An Interpolation Algorithm}

Applying Algorithm G to the set of generators of $I_{v,m,l}$ in Corollary \ref{propskqq}, we obtain an interpolation algorithm for the list decoding of Reed-Solomon codes.

\textbf{Algorithm I.} Given input $v=(v_1,v_2,\dots,v_n)$ and parameters $m$ and $l$, this algorithm finds the minimal polynomial of $I_{v,m,l}$ with respect to monomial order $>_{k-1}$. Let $g_i=\sum_{j=0}^la_{ij}y^j$ for $0\le i\le l$ during the execution of the algorithm.

\begin{enumerate}
\item[I1.] Compute $h_v=\sum_{i=1}^nv_ih_i$.
\item[I2.] For $0\le i\le m$, set
\[
	g_i\leftarrow(y-h_v)^i\prod_{j=1}^n(x-\ga_j)^{m-i}
\] 
and for $m<i\le l$, set 
\[
	g_i\leftarrow y^{i-m}(y-h_v)^m. 
\]
Set $r\leftarrow 0$.
\item[I3.] Increase $r$ by $1$. If $r\le l$, then proceed; otherwise go to step I7.
\item[I4.] Find $s=\ydeg(\LT(g_r))$. If $s=r$, then go to step I3.
\item[I5.] Set $d\leftarrow\deg(a_{rs})-\deg(a_{ss})$ and  $c\leftarrow\LC(a_{rs})\LC(a_{ss})^{-1}$.
\item[I6.] If $d\ge 0$, then set 
\[
	g_r\leftarrow g_r-cx^dg_s.
\]
If $d<0$, then set, storing $g_s$ in a temporary variable,
\[
	g_s\leftarrow g_r,\quad g_r\leftarrow x^{-d}g_r-cg_s.
\]
Go back to step I4.
\item[I7.] Let $Q$ be the $g_i$ with the smallest leading term. Output $Q$ and the algorithm terminates.  
\end{enumerate}

\begin{exmp}\normalfont 
Let $\F_7=\set{0,1,2,\dots,6}$ be the finite field with $7$ elements. Let $n=6$ and $k=3$. Choose $\ga_i=i$ for $1\le i\le 6$.  We use $\RS(6,3)$ over $\F_7$ as an example. Suppose that $v=(6,2,4,4,4,2)$ is the received vector. 

Let $m=2$, and consider $I_{v,2}$. Let $Q$ be the minimal polynomial of $I_{v,2}$ with respect to $>_2$. We set $l=3\ge\ydeg(Q)$. For our $v$, we have
\[
	h_v=x^4+5x^3+4x^2+4x+6,\quad\eta=\prod_{j=1}^n(x-\ga_j)=x^6-1.
\]
Therefore
\[
\begin{split}
	I_{v,2,3}&=\langle \eta^2,(y-h_v)\eta,(y-h_v)^2,y(y-h_v)^2\rangle\\
		&=\langle\eta^2,\eta y-\eta h_v,y^2-2h_vy+h_v^2,y^3-2h_vy^2+h_v^2y\rangle.
\end{split}
\]
Let $g_0=\eta^2$, $g_1=\eta y-\eta h_v$, and so on. Note that $\ydeg(g_i)=i$ for $0\le i\le 3$. 

We demonstrate Algorithm I by finding the minimal polynomial of $I_{v,2,3}$. In the following, polynomials in $\F[x]$ are parenthesized with only leading terms shown. After steps I1 and I2, we have
\[
	\left\{\begin{array}{crrrrrrrrrrrrr}
	g_0=&&&&&&&(x^{12}+\cdots)\\
	g_1=&&&&&(x^6+\cdots)y&+&(6x^{10}+\cdots)\\
	g_2=&&&y^2&+&(5x^4+\cdots)y&+&(x^8+\cdots)\\
	g_3=&y^3&+&(5x^4+\cdots)y^2&+&(x^8+\cdots)y \\
	\end{array}\right.
\]
After step I3, $r=1$. In step I4, we find $s=\ydeg(\LT(g_1))=0$. Since $s\neq r$, we go to step I5. Then $d=-2$ and $c=6$. So in step I6, $g_0$ and $g_1$ is replaced with $g_1$ and $x^2g_1-6g_0$, respectively. Then we have
\[
	\left\{\begin{array}{crrrrrrrrrrrrr}
	g_0=&&&&&(x^6+\cdots)y&+&(6x^{10}+\cdots)\\
	g_1=&&&&&(x^8+\cdots)y&+&(2x^{11}+\cdots)\\
	g_2=&&&y^2&+&(5x^4+\cdots)y&+&(x^8+\cdots)\\
	g_3=&y^3&+&(5x^4+\cdots)y^2&+&(x^8+\cdots)y \\
	\end{array}\right.
\]
After one more update like this, we have
\[
	\left\{\begin{array}{crrrrrrrrrrrrr}
	g_0=&&&&&(x^6+\cdots)y&+&(6x^{10}+\cdots)\\
	g_1=&&&&&(x^8+\cdots)y&+&(2x^9+\cdots)\\
	g_2=&&&y^2&+&(5x^4+\cdots)y&+&(x^8+\cdots)\\
	g_3=&y^3&+&(5x^4+\cdots)y^2&+&(x^8+\cdots)y \\
	\end{array}\right.
\]
This time we find $s=\ydeg(\LT(g_1))=1$. Since $s=r$, we go to step I3, and increase $r$ by one. In step I4, we find $s=\ydeg(\LT(g_2))=0$. Since $s\neq r$, we go to step I5. Then $d=-2$ and $c=6^{-1}=6$. So in step I6, $g_0$ and $g_2$ is replaced with $g_2$ and $x^2g_2-6g_0$, respectively. Then we have
\[
	\left\{\begin{array}{crrrrrrrrrrrrr}
	g_0=&&&y^2&+&(5x^4+\cdots)y&+&(x^8+\cdots)\\
	g_1=&&&&&(x^8+\cdots)y&+&(2x^9+\cdots)\\
	g_2=&&&x^2y^2&+&(6x^6+\cdots)y&+&(5x^9+\cdots)\\
	g_3=&y^3&+&(5x^4+\cdots)y^2&+&(x^8+\cdots)y \\
	\end{array}\right.
\]
The algorithm continues updating in the same way. After the final update, we have
\[
	\left\{\begin{array}{crrrrrrrrrrrrr}
	g_0=&&&y^2&+&(5x^4+\cdots)y&+&(x^8+\cdots)\\
	g_1=&&&(x^2+\cdots)y^2&+&(6x^6+\cdots)y&+&(5x^7+\cdots)\\
	g_2=&y^3&+&(6x^3+\cdots)y^2&+&(3x^5+\cdots)y&+&(4x^7+\cdots)\\
	g_3=&xy^3&+&(4x^3+\cdots)y^2&+&(3x^5+\cdots)y&+&(6x^6+\cdots) \\
	\end{array}\right.
\]
This set $\set{g_0,g_1,g_2,g_3}$ is a Gr\"obner basis of $I_{v,2,3}$. The minimal polynomial is $g_2$. So the algorithm terminates with output
\begin{multline*}
	Q=y^3+(6x^3+4x+5)y^2+(3x^5+6x^4+4x^3+6x^2+6x+2)y\\+4x^7+4x^6+3x^5+3x^4+4x^3+2x^2+x+6.
\end{multline*}
Since $Q$ has factorization
\[
  (y+x^2+5x+2)(y+3x^2+4x+6)(x^3+6y+4x^2+5x+3),
\]
a root-finding algorithm will output $6x^2+2x+5$ and $4x^2+3x+1$ with degree less than $3$, each of which yields a codeword $c$ with $\wt(v-c)\le 2$.
\end{exmp}

\section{Complexity of the Interpolation Algorithm}

We give an upper bound on the number of multiplication operations in the field $\F$ required during the execution of Algorithm I. We assume that the operation of polynomial multiplication is done in the straightforward method such that a multiplication of two polynomials of degree $a$ and $b$ requires $(a+1)(b+1)$ multiplication operations over $\F$.

Step I1 requires at most 
\[
	n^2+\sum_{i=2}^mn((i-1)(n-1)+1)=O(n^2m^2)
\]
multiplication operations. Step I2 requires at most
\[
	\sum_{i=0}^m\sum_{j=0}^i((i-j)(n-1)+1)((m-i)n+1)=O(n^2m^4)
\]
multiplication operations. To analyze the iterative steps I3--I6, fix $0\le r\le m$. Observe that at the start of the updating for $g_r$, the leading term of $g_r$ is in $a_{r0}$, and
\[
	\deg(a_{r0})-\deg_{k-1}(a_{rr}y^r)=\deg((-h_v)^r\eta^{m-r})-\deg_{k-1}(\eta^{m-r}y^r)\le(n-k)r.
\]
Then Corollary \ref{fdhwkd} implies that at most $(n-k)r^2$ updates take place for $r$. Hence the total number of updates for all $0\le r\le m$ is
\[
	\sum_{r=0}^m(n-k)r^2.
\]
For each update, step I6 requires at most
\[
	\sum_{j=0}^r(mn-j(k-1)+1)
\]
multiplication operations, because it always holds that $\deg_{k-1}(g_i)\le mn$ for $0\le i\le m$. For $m<r\le l$, we can do a similar analysis. To summarize, steps I3--I6 take totally at most 
\[
\begin{split}
	&\sum_{r=0}^m\sum_{j=0}^r(n-k)r^2(mn-j(k-1)+1)\\
	&\quad+\sum_{r=m+1}^l\sum_{j=0}^r(n-k)mr(mn+(r-m-j)(k-1)+1)\\
	&=O(n^2ml^4)
\end{split}
\]
multiplication operations. As $l$ can be set to $O(mn^{1/2}k^{-1/2})$ by \eqref{dmxcxn}, we conclude that an execution of Algorithm I takes $O(n^4k^{-2}m^5)$ multiplication operations over $\F$.

\section{A Special Case}

Let us consider Guruswami and Sudan's list decoding for the case $m=l=1$. In this case, our interpolation algorithm becomes simplest and a root-finding algorithm is not necessary for decoding. Thus we obtain a simple decoding algorithm of Reed-Solomon codes. 

Let $I_v=I_{v,1,1}$. We begin with considering the minimal polynomial of $I_v$ with respect to $>_{k-1}$. Let $ay+b$ be the minimal polynomial with $a,b\in\F[x]$. We want to have an upper bound on the $(1,k-1)$-weighted degree of $ay+b$. Proposition \ref{propsjqx} implies that the monomials occurring in $ay+b$ belong to the first $n+1$ monomials of $\F[x,y]_1$ in the order $>_{k-1}$. So we consider the following table of monomials of $\F[x,y]_1$ ordered in $>_{k-1}$
\[
\begin{array}{rrrrr|rrrrrrr}
	&&&&&y&xy&\cdots&x^{k-2}y&x^{k-1}y&x^ky&\cdots\\
	1&x&x^2&\cdots&x^{k-2}&x^{k-1}&x^k&\cdots&x^{2k-3}&x^{2k-2}&x^{2k-1}&\cdots
\end{array}
\]
where the ordering is from left to right and from bottom to top. Consider the first $n+1$ monomials in the table. Let us index only the columns of length two so that the column containing $x^{k-1}$ has index $0$. Let $C$ be the index of the column in which $(n+1)$-th monomial lies. Then $C$ is the smallest integer satisfying
\[
	k-1+2(C+1)\ge n+1,
\]
namely $C=\lceil(n-k)/2\rceil$. It follows that every monomial occurring in $ay+b$ has $(1,k-1)$-weighted degree $\le k-1+C$. We conclude that $\deg_{k-1}(ay+b)\le k-1+\lceil(n-k)/2\rceil$.

Proposition \ref{propdjqq} allows us to exactly determine the form of the minimal polynomial of $I_v$ with respect to $>_{k-1}$.

\begin{prop}\label{propsjqq}
Let $\tau=\lfloor(n-k)/2\rfloor$. There is at most one codeword $c$ satisfying $\wt(v-c)\le\tau$. Suppose that there is such a codeword $c$. Let $e=v-c$, and 
\[
	f_e=\prod_{e_i\not=0}(x-\ga_i).
\]
Then $f_e(y-h_c)$ is the minimal polynomial of $I_v$ with respect to $>_{k-1}$.
\end{prop}

\begin{pf}
Let $ay+b$ be the minimal polynomial of $I_v$ with respect to $>_{k-1}$. Let $w=\deg_{k-1}(ay+b)$. Since $w\le k-1+\lceil(n-k)/2\rceil$, we have
\[
	n-w-1\ge n-k-\lceil(n-k)/2\rceil=\tau.
\]
Then Proposition \ref{propdjqq} says that every codeword $c$ satisfying $\wt(v-c)\le\tau$ yields a root $h_c$ of $ay+b$. Since $ay+b$ can have at most one root, it follows that there is at most one codeword $c$ satisfying $\wt(v-c)\le\tau$.

Suppose that $c$ is such a codeword. Then $ay+b=a(y-h_c)$. Let $e=v-c$. Since $a(y-h_c)\in I_v$, for each $1\le i\le n$,  
\[
	0=a(\ga_i)(v_i-h_c(\ga_i))=a(\ga_i)e_i.
\]
When $e_i\neq 0$, we must have $a(\ga_i)=0$. Thus $f_e$ divides $a$. Since $f_e(y-h_c)\in I_v$, the minimality of $ay+b$ implies that $ay+b=f_e(y-h_c)$.
\end{pf}

We now assume that there occurred no more than $\tau=\lfloor(n-k)/2\rfloor$ errors to the sent codeword. Then Proposition \ref{propsjqq} says that the sent codeword $c$ is the unique codeword satisfying $\wt(v-c)\le\tau$, and the message polynomial $h_c$ is obtained by one division from the minimal polynomial of $I_v$. On the other hand, Algorithm I is substantially simplified when it is applied to $I_v=I_{v,1,1}$. Hence we have the following

\bigskip
\textbf{Decoding Algorithm D.} Given the received vector $v=(v_1,v_2,\dots,v_n)$, this algorithm finds the message polynomial $h_c$. The polynomials $\eta=\prod_{j=1}^n(x-\ga_j)$ and $h_i$ as in \eqref{equcss} for $1\le i\le n$ are precomputed.

\begin{enumerate}
\item[D1.] Compute $-h_v=-\sum_{i=1}^nv_ih_i$.
\item[D2.] Set
\[
\begin{aligned}
	A\leftarrow 0,\quad B\leftarrow \eta,\quad
	C\leftarrow 1,\quad D\leftarrow -h_v.
\end{aligned}
\]
\item[D3.] If $\deg(C)+k-1\ge\deg(D)$, then go to step D6.
\item[D4.] Set $d\leftarrow\deg(D)-\deg(B)$ and  $c\leftarrow\LC(D)\LC(B)^{-1}$.
\item[D5.] If $d\ge 0$, then set 
\[
	C\leftarrow C-cx^dA,\quad
	D\leftarrow D-cx^dB.
\]
If $d<0$, then set, storing $A$ and $B$ in temporary variables,
\[
\begin{aligned}
	A\leftarrow C,
		&\quad B\leftarrow D,\quad
	C\leftarrow x^{-d}C-cA,
		&\quad D\leftarrow x^{-d}D-cB.
\end{aligned}
\]
Go back to step D3.
\item[D6.] Set $h\leftarrow -D/C$. Output $h$ and the algorithm terminates.  
\end{enumerate}

Recall that generalized Reed-Solomon codes are defined as a simple twist of Reed-Solomon codes. Hence it is straightforward  to modify our decoding algorithm to work for generalized Reed-Solomon codes as well. Then the modified algorithm decodes alternant codes up to half of the designed distance, as alternant codes are defined as subfield subcodes of generalized Reed-Solomon codes. For example the modified algorithm decodes BCH codes up to half of the designed distance. We leave the details to the reader.

\section{Conclusion}

We focused on the interpolation problem in Guruswami and Sudan's list decoding of Reed-Solomon codes. Though we are well aware of the important extension of their idea for soft-decision decoding, we restricted our attention to hard-decision decoding where multiplicities are assigned uniformly. Here we just note that our results can be easily extended for soft-decision decoding of Reed-Solomon codes by finding a suitable set of generators of the ideal of interpolation polynomials for arbitrary points with arbitrary multiplicities. See \citet{kwankyu2006b} for an extension in this direction.

For the problem of computing a Gr\"obner basis of the vanishing ideal of points with multiplicities on the plane, common wisdom is to use Buchberger's algorithm or the Marinari-M\"oller-Mora algorithm in \citet{marinari1993}. However, for the application to decoding, either algorithm needs to be optimized exploiting the particular need of finding the $Q$-polynomial of the interpolation ideal, rather than the whole Gr\"obner basis, with respect to the particular weighted monomial order. Here we presented such an optimized version of Buchberger's algorithm, though our presentation is self-contained and an explicit complexity analysis is given.

One may notice some similarities between our algorithm computing a Gr\"obner basis of a module over a univariate polynomial ring and the algorithm of \cite{aleknovich2005} computing a reduced basis of a lattice over a univariate polynomial ring. Moreover, to compute the minimal polynomial of the interpolation ideal, both algorithms rely on a set of generators of the ideal. However, working with the module induced from the interpolation ideal, our interpolation algorithm computes the minimal polynomial of the ideal more directly and systematically than Alekhnovich's algorithm. We remark that our module-theoretic approach was inspired by the illuminating work of \cite{fitzpatrick1995}.

\begin{ack}
We thank the referees for the useful suggestions that greatly improved the exposition. 

The first author was supported by the Korea Research Foundation Grant funded by Korea Government (MOEHRD, Basic Research Promotion Fund) (KRF-2005-214-C00009).
\end{ack}


\end{document}